\input amstex
\documentstyle{amsppt}
\NoBlackBoxes

  {}
  {}

\redefine\sp{Sprind\v zuk}

\define\qu{quasiunipotent}
\define\hd{Hausdorff dimension}
\define\hs{homogeneous space}
\define\ggm{G/\Gamma}
\define\df{\overset\text{def}\to=}
\define\br{\Bbb R}
\define\bn{\Bbb N}
\define\bz{\Bbb Z}
\define\bq{\Bbb Q}
\define\ba{badly approximable}

\define\di{Diophantine}

\define\bc{\Bbb C}

\define\da{Diophantine approximation}

\define\ehs{expanding horospherical subgroup}

\define\va{\bold b}
\define\ve{\bold e}
\define\vf{\bold f}
\define\vx{\bold x}
\define\vv{\bold v}
\define\vy{\bold y}

\define\vt{\bold t}
\define\vw{\bold w}
\define\vp{\bold p}
\define\vq{\bold q}

\define\mn{{m+n}}

\define\mr{M_{m\times n}(\br)}
\define\amr{$Y\in\mr$}
\define\la{L_{ Y}}
\define\ly{L_{\vy}}
\redefine\aa{\langle Y,\va\rangle}

\define\mtr{\tilde M_{m \times n}(\br)}

\define\amtr{$\aa\in \mtr$}

\define\nz{\smallsetminus \{0\}}
\define\vre{\varepsilon}
\define\ssm{\smallsetminus}
\def\Ad{\operatorname{Ad}}
\def\diag{\operatorname{diag}}
\def\const{\operatorname{const}}
\def\Id{\operatorname{Id}}
\define\ph{partially hyperbolic}
\define\cag{$(C,\alpha)$-good}


\topmatter
\title Some applications of homogeneous dynamics to number theory \\
\endtitle 

\author {Dmitry Kleinbock} 
\endauthor


    \address{ Dmitry  Kleinbock,  Department of
Mathematics, Brandeis University, Waltham, MA 02454-9110}
  \endaddress

\email kleinboc\@brandeis.edu \endemail
\thanks The author was supported in part by NSF
Grants DMS-9704489 and DMS-0072565.\endthanks

%


\endtopmatter


\document

This survey paper is not a complete reference guide to
number-theoretical 
applications of ergodic theory. Instead, the plan is to consider an
approach to a class of problems involving \di\ properties of
$n$-tuples of real numbers, namely, describe a specific dynamical
system which is naturally connected with these problems. 

\head{1. A glimpse at  \da}\endhead


For motivation, let us  start by looking at two  (vaguely
defined) \di\ problems: 

\proclaim{Problem 1}
Given a nondegenerate 
indefinite quadratic form  of
signature $(m,n)$, study 
the set of 
its values at integer points.  
\endproclaim

Here is a precise statement along these lines, conjectured in 1929 by
Oppenheim \cite{Op1} and proved in 1986 by Margulis \cite{Ma3}:

\proclaim{Theorem 1.1} Let $B$ be a real  nondegenerate indefinite
quadratic form  of
signature $(m,n)$, $k = m+n > 2$. Then either $B$ is proportional to a
rational form, or $\inf_{\vx\in\bz^{k}\nz}|B(\vx)| = 0$. 
\endproclaim

One possible approach to the problem is to  write 
$
B(\vx) = \lambda
S_{m,n}(g\vx)\,,
$
 where $\lambda\in\br$, $g\in SL_k(\br)$ and 
$$
S_{m,n}(x_1,\dots,x_k) = x_1^2 + \dots + x_m^2
- x_{m+1}^2 - \dots - x_k^2
$$
(a linear unimodular change of variables).  Then the problem reduces
to studying values of the standard form $S_{m,n}$ of signature $(m,n)$
applied to the collection of 
vectors of the form $\{g\vx \mid \vx
\in \bz^k\}$. And the dynamical approach consists of studying the
action of the stabilizer of the form $S_{m,n}$ on such collections.

\proclaim{Problem 2} 
Given  $m$ vectors $\vy_1,\dots,\vy_m\in\br^n$ (viewed as linear forms
$\vx\mapsto \vy_i\cdot\vx$,  $\vx\in\br^n$)
 how small (simultaneously) can be the values of  
$
|\vy_i\cdot\vq + p_i|, \ p_i\in\bz, 
$
when $\vq = (q_1,\dots,q_n)\in\bz^n$ is far from $0$? 
\endproclaim

Let us also illustrate this by a conjecture, this time still
open. Here we specialize to the case of just one linear form
given by $\vy\in\br^n$. The following is
known as Littlewood's (1930) Conjecture:

\proclaim{Conjecture 1.2} For every $\vy\in\br^n$, $n \ge 2$, one has 
$$
\inf_{\vq\in\bz^n\nz,\,p\in\bz}|\vy\cdot\vq +
p|\cdot\Pi_{\sssize +}(\vq) = 0\,,\tag 1.1
$$
where $\Pi_{\sssize +}(\vq)$ is defined to be equal to $\prod_{i =
1}^n \max(|q_i|, 1)$ or, 
equivalently, $\prod_{q_i \ne 0} |q_i|$.
\endproclaim

To approach Problem 2, one can put together 
$$
\vy_1\cdot\vq + p_1, \dots,\vy_m\cdot\vq + p_m \quad \text{and}\quad 
q_1,\dots,q_n,
$$
 and consider the {\it lattice}
$$
\left\{\left.\left(\matrix Y\vq + \vp \\  \vq \endmatrix
\right)\right| \vp \in\bz^m,\,\vq \in\bz^n\right\} = L_{
Y}\bz^{\mn}\,,\tag 1.2 
$$
where 
$$
L_{ Y} \df \left(\matrix
I_m & Y \\
 0 & I_n
\endmatrix \right)\tag 1.3
$$ and $Y$ is the matrix with rows $\vy_1^{\sssize T},\dots,\vy_m^{\sssize T}$. In this 
case, the orbit of the lattice (1.2) under a certain group action
provides a way to study Diophantine properties of 
$\vy_1,\dots,\vy_m$. 

\medskip 

In both cases, we use the initial data of a number-theoretic problem
to construct a  lattice in a Euclidean space, and then work with
the collection of all such objects (lattices). 
Our goal is to describe several principles
responsible for a particular class of applications of flows in the
space of lattices to number theory. For more details and a broader
picture the reader is referred to a number of extensive reviews of
homogeneous actions and interactions with number theory which have
appeared during the last 10 years, such as: ICM talks of Margulis
\cite{Ma5}, 
Ratner \cite{Ra3}, Dani \cite{D5} and Eskin \cite{E}, books
\cite{St2} and \cite{BMa}, and survey 
papers \cite{D6, D7, KSS, Ma7, St1}. 

The structure of this paper is as follows: in the next section we
collect all basic facts about the space of lattices, and discuss
a ``lattice''  approach to  studying values of quadratic forms at
integer points (Problem 1 and Theorem 1.1 in particular). Then in
\S\S\ 3 and  4 we take several sub-problems of Problem 2 and
describe recent 
results obtained by means of homogeneous dynamics. The last section is
devoted to Conjecture 1.2 and related issues, that is, so called
``multiplicative \da''.

\head{2. The space of lattices}\endhead

\subhead  {Phase space}
\endsubhead 
Fix $k \in
\bn$ and consider
$$
\Omega \df \text{the set of unimodular lattices in } \br^k
$$
(discrete subgroups with covolume $1$).
That is, any lattice $\Lambda\in\Omega$ is equal to $\bz\vv_1 \oplus \dots
\oplus \bz\vv_k$, where the set $\{\vv_1,\dots,\vv_k\}$ (called a {\sl
generating set\/} of the 
lattice) is  linearly independent, and the volume of the
parallelepiped spanned by $\vv_1,\dots,\vv_k$ is equal to $1$.

An element of $\Omega$ which is easy to distinguish is $\bz^k$ (the
standard lattice). In fact, any $\Lambda\in \Omega$ is equal to
$g\bz^k$ for some $g\in G \df SL_k(\br)$. That is, $G$ acts
transitively on $\Omega$, and, further, $\Gamma \df SL_k(\bz)$ is the
stabilizer of $\bz^k$. In other words, $\Omega$ is isomorphic to the
\hs\ $G/\Gamma$. 

\subhead  {Measure}
\endsubhead 
One can consider a Haar measure on $G$ (both left and right invariant)
and the corresponding left-invariant measure on $\Omega$. It is well
known that  the
resulting measure happens to be finite. We
denote by $\mu$ the normalized Haar measure on $\Omega$.

\subhead  {Topology}
\endsubhead Two lattices are said to be close if generating sets which are
close to each other can be chosen for them.  
This defines a topology on $\Omega$ which coincides with the
quotient topology on $G/\Gamma$.  An important feature is that $\Omega$ is
not compact (in other words, $\Gamma$ is a {\sl non-uniform\/} lattice
in $G$). More precisely, one has

 \proclaim{Theorem 2.1 {\rm (Mahler's Compactness
Criterion, see \cite{R})}} A subset $K$ of $\Omega$ is bounded iff
there exists $\vre > 0$ such that for any $\Lambda \in K$ one has
$\inf_{\vx\in\Lambda\nz}\|\vx\| \ge \vre$. In other words, define 
$$
\Omega_\vre \df\big\{\Lambda \in \Omega \bigm|  \|\vx\| < \vre\text{ for some
}\vx\in\Lambda\nz\big\}\,;
$$
then $\Omega\ssm\Omega_\vre$ is compact. 
\endproclaim

\subhead  {Action}
\endsubhead $\Omega$ is a topological $G$-space, with the
(continuous) left action defined by 
$$
g\Lambda = \{g\vx\mid \vx\in\Lambda\} \quad\text{or} \quad g(h\Gamma) = (gh)\Gamma\,.
$$
One can consider the
 action of various subgroups (one- or multi-parameter) or
subsets of 
$G$. Thus one gets an interesting class of dynamical systems. 
Several important features of these systems are worth mentioning.  

First,
the geometry of the phase space is ``uniform'': a small enough neighborhood of
every point of $\Omega$ is isometric to  a
neighborhood of identity in  $G$. In other words, many geometric
constructions can be reduced to algebraic manipulations in $G$. 

Second,  the very rich representation theory of $G$ can be heavily
used. Namely,
the $G$-action on $\Omega$ can be studied via the regular representation of
$G$ on $L^2(\Omega)$

The two features above in fact apply for all {\sl homogeneous
actions\/}, that is, actions of subgroups of a Lie group $G$ on the
quotient space $\ggm$ where $\Gamma$ is a lattice in $G$. There are
also important features specifically for the space $\Omega =
SL_k(\br)/SL_k(\bz)$: namely, combinatorial structure of the space of
lattices, as well as intuition coming from the theory of \da. 
\vskip .1in 

In what follows we will focus our attention on the space $\Omega$, but
most of the results will be valid in much bigger generality of
homogeneous actions, which
will be indicated. The reader
is referred to \cite{AGH,
Ma6, R, St2, Z} for general
facts about Lie groups,
discrete subgroups and \hs s. 

\subhead  {Classification of actions} \endsubhead Let $G$ be a
Lie group and  
$\Gamma$ a discrete subgroup. 
Since $g(h\Lambda ) = (ghg^{-1})g\Lambda $ for every $\Lambda \in\ggm$
and $g,h\in G$, local properties of the
$g$-action are determined by the differential of the conjugation map,
$\Ad_g(x)=\frac{d \big(g\exp(tx)g^{-1}\big)}{dt}|_{t=0}$ (here $x$
belongs to the Lie algebra of $G$).
An element $g\in G$ is said to be:
{\sl unipotent} if  $(\Ad_g-\Id)^j=0$ for some $j\in\bn$ (equivalently,  all eigenvalues of $\Ad_g$ are equal to $1$);
{\sl quasi-unipotent} if  all eigenvalues of $\Ad_g$ are of absolute value $1$;
{\sl \ph} if it is not quasi-unipotent.

Given $g\in G$, define 
$$
H^{\pm}(g) = \{h\in G\mid g^{-l}hg^{l}\to e \text{ as } l\to \pm\infty\}
$$
({\sl  expanding and contracting horospherical subgroups}). 
Then $G$ is locally a direct product of $
H^{-}(g)$, $
H^{+}(g)$ and another subgroup $H^0(g)$,  and $g$ is \qu\ iff $H^0(g) =
G$ (that is, $
H^{-}(g)$ and $
H^{+}(g)$ are trivial). 
Furthermore, for any $\Lambda \in\ggm$ the orbits $
H^{-}(g)\Lambda $,  $
H^{+}(g)\Lambda $ and $H^0(g)\Lambda $ are leaves of {\sl stable,
unstable and neutral\/} foliations on $\ggm$. 

We now specialize to the case $G = SL_k(\br)$ and $\Gamma = SL_k(\bz)$.

\example{Example}
The simplest case is when $k = 2$: then $\Omega = \ggm$ is isomorphic to  the unit tangent bundle
to the surface $\Bbb H^2/SL_2(\bz)$. 
The
geodesic flow  is then given by   the action of
$\pmatrix e^t & 0  \\ 0 & e^{-t}  \endpmatrix$, and the
horocycle flow  -- by    the action of
$\pmatrix 1 & t  \\ 0 & 1  \endpmatrix$ (the simplest example of
a unipotent flow). 

\endexample

\example{More examples}
Suppose that $g\in G$ is diagonalizable over
$\br$, and take a basis of $\br^k$ in which $g =
\diag(\underbrace{\lambda_1,\dots,\lambda_1}_{\text{$i_1$
times}},\dots\dots,\underbrace{\lambda_l,\dots,\lambda_l)}_{\text{$i_l$
times}},\quad\lambda_1 > \dots > \lambda_l$. 
Then $
H^{\sssize -}(g)$ 
and  $
H^{\sssize +}(g)$ 
are subgroups     
of lower- 
and
upper- 
 triangular 
 groups. An
important special case occurs when $g$ as above comes from ``the most
singular'' 
direction in a Weyl chamber of the Lie algebra of
$G$; that is, when it has only two distinct eigenvalues. In this case
one can write $k = \mn$ and consider a one-parameter subgroup of
$G = SL_k(\br)$  given by 
$$
g_t = \text{\rm
diag}(e^{t/m},\dots,e^{t/m},e^{-t/n},\dots,e^{-t/n})\,.\tag 2.1
$$
 Then the \ehs\ of $G$ relative to $g_1$ is exactly $\{\la \mid
Y\in\mr\}$, where $\la$ is as defined in (1.3). 
\endexample

\subhead  {Ergodic properties} 
\endsubhead
Here the main tool is the representation theory of semisimple Lie
groups. 
By a theorem 
of Moore \cite{Mo1}, the action of any noncompact closed subgroup
of  $G$  on $\Omega = \ggm$  
is ergodic and, moreover, 
mixing; in other words, matrix coefficients $(g\varphi, \psi)$ of
square-integrable 
functions on $\Omega$ with mean value zero tend to $0$ as $g\to\infty$
in $G$. (Here $(\cdot,\cdot)$ stands for the inner product in $L^2(\Omega)$.) In fact for smooth functions this decay is exponential, as
shown in the following 
 
 \proclaim{Theorem 2.2 (Decay of correlations)}  There exists
$\beta > 0$ such that  
 for any two functions $\varphi, \psi
\in C^\infty_{comp}(\ggm)$   with $\int
\varphi
\,d\mu 
=
\int
\psi
\,d\mu 
= 0$ and any $g\in G$ one has
$$
\big|\int
(g\varphi
\cdot\psi)\,d\mu 
\big|\le \const(\varphi, \psi)e^{-\beta\|{g}\|}\,.
$$
In particular, if $g_t$ is \ph, then
$$
\big|\int
(g_t\varphi
\cdot\psi)\,d\mu 
\big|\le \const(\varphi, \psi,g_t)e^{-\gamma t}\,.\tag 2.2
$$
\endproclaim 

See \cite{Mo2, Ra1} for $k = 2$, \cite{KS} for $k >
2$. 
\vskip .1in 

The following result can be derived from the mixing property of
partially hyperbolic actions on $\ggm$:

\proclaim{Theorem 2.3 {\rm (Uniform distribution of unstable leaves,
\cite{KM1})}} 
Let $g_t$ be a
partially hyperbolic one-parameter subgroup of $G$, $H = H^{\sssize
+}(g_1)$,  $\nu$ a Haar
measure on $H$. Then 
for any open subset $V$
of $H$, any $\varphi\in C^\infty_{comp}(\ggm)$ and any compact
subset $Q$ of $\ggm$, the average of $\varphi$ over the $g_t$-image
of $V\Lambda $, $\Lambda \in\ggm$, tends to the integral of  $\varphi$ as $t\to
\infty$ uniformly (in $\Lambda $) on compact subsets of $\ggm$; that
is, 
$$
 \frac1{\nu(g_tVg_{-t})} \int_{g_tVg_{-t}}\varphi(hg_t\Lambda )\,d\nu(h)
\to \int_{\ggm}\varphi\,d\mu\,.\tag 2.3
$$
\endproclaim

\remark{Remarks} Moore's theorem (that is, a criterion for
mixing of subgroup actions) was proved under the
assumption that
$$
\aligned
\text{ $G$ is a connected semisimple Lie
group with finite center}\\ \text{ and no compact factors, 
and $\Gamma$ is an
irreducible lattice in }&G\,,\endaligned \tag 2.4
$$
(A lattice is irreducible if it is not, up
to commensurability, a product of lattices in simple factors of $G$.) See also \cite{BM, Ma5} for more general
ergodicity and mixing criteria. 
Theorems 2.2 was proved in \cite{KS} assuming (2.4) 
and in addition  that 
$$
\text{all simple factors of $G$ have property (T).}\tag 2.5
$$
One also knows, see \cite{Bek, \rm Lemma 3}, that Theorem 2.2 holds when the group
$G$ is simple. In
\cite{KM3} it was shown that one can remove condition (2.5) but
instead assume 
that $\Gamma\subset G$ is a non-uniform lattice. 

Theorem 2.3 is also
proven in the generality of the assumption (2.4) (in fact, one only
needs mixing of the $g_t$-action), and as long as Theorem 2.2 holds,
the convergence in (2.3) is 
exponential in $t$. See \cite{KM1, K4} for more details and
generalizations, and \cite{Ma7, \rm Remark 3.10} for references to other
related results and methods. 
\endremark

\subhead  {Recurrence of unipotent trajectories} 
\endsubhead
In this subsection we are back to the case $G = SL_k(\br)$ and $\Gamma
= SL_k(\bz)$. It is an  elementary geometric observation that
horocyclic trajectories on 
$SL_2(\br)/SL_2(\bz)$ do not run off to infinity. 
It is much harder to prove that the same holds for any unipotent flow on
$\Omega = SL_k(\br)/SL_k(\bz)$ for $k \ge 3$ {\rm \cite{Ma2, D3}}. The
theorem below, due to Dani (1985), is a 
quantitative strengthening; it shows that for any unipotent orbit one
can find a compact subset of $\Omega$ such that the density of time
that the orbit spends in this set is as close to $1$ as one wishes:

 \proclaim{Theorem 2.4 \cite{D3}}  For any $\Lambda\in\Omega$
and any $\delta > 0$ 
there exists $\vre > 0$ such that  for any unipotent subgroup
$\{u_x|x\in\br\}$ of $G$ and any $T > 0$ one has
$$
\big|\{x\in[0,T]\mid u_x\Lambda\in\Omega_\vre\}\big| \le \delta T \,.
$$
\endproclaim
 \vskip .1in 

The proof is based on the combinatorial structure of the space of
lattices. We will obtain the theorem above as a corollary of a more
general fact in \S 4. 

\subhead  {Orbit closures of unipotent flows} 
\endsubhead
It has been proved by Hedlund that  any orbit of the horocycle flow on 
$SL_2(\br)/SL_2(\bz)$ is either periodic or dense. A
far-reaching generalization has been conjectured by Raghunathan and
proved in full generality by 
Ratner. In particular, one has the following

 \proclaim{Theorem 2.5} {\rm \cite{Ra2}} Let $G$ be a connected Lie
group, $\Gamma$ a lattice in $G$, and let $U$ be a subgroup of
$G$ generated by unipotent one-parameter subgroups. Then for any
$x\in\ggm$ there exists a closed subgroup 
$L$ containing $U$ such that the closure of the orbit
$Ux$
coincides with $Lx$ and there is an $L$-invariant probability
measure supported on $Lx$. 
\endproclaim
 
We emphasize that Theorem 2.5 allows one to understand all (not just
almost all) orbits. It is this feature which is responsible for 
applications of this theorem to number theory. Let us illustrate it by
sketching the reduction of the Oppenheim Conjecture (Theorem 1.1) to
Ratner's theorem.

\proclaim{Corollary 2.6} Let $
S(x_1,x_2,x_3) = 2x_1x_3 - x_2^2\,,
$ and 
$$
H_S = \{h\in SL_3(\br)\mid S(h\vx) = S(\vx)\
\forall\,\vx\in\br^3\}\cong 
SO(2,1) 
$$ 
(the stabilizer of $S$). Then any relatively compact orbit
$H_S\Lambda$, $\Lambda$ a lattice in 
$\br^3$, is compact.
\endproclaim

\demo{Proof} $H_S$ is generated by its unipotent
one-parameter subgroups, namely
$$
u(t) =\pmatrix 1 &t&t^2/2\\ 0 &1 &t\\ 0 &0 & 1\endpmatrix \text{
and }u^{T}(t) =\pmatrix 1 &0&0\\ t &1 &0\\ t^2/2 &t & 1\endpmatrix\,, 
$$
and 
there are no intermediate subgroups between  $H_S$ and $SL_3(\br)$.
Hence by Theorem 2.5 any $H_S$-orbit is either closed or
dense. \qed\enddemo 

The following is crucial for the deduction of the Oppenheim conjecture
from the above corollary (implicitly stated in \cite{CS} and later
observed by Raghunathan): 

\proclaim{Lemma 2.7} Let $B$ be a real  nondegenerate indefinite
quadratic form in $3$ 
variables. Write  $
B(\vx) = \lambda
S(g\vx)$ for some $g\in SL_3(\br)$. Then  the orbit
$H_S g\bz^3$ is relatively compact if and only if
$$
|B(\vx)| \ge \vre \text{ for some $\vre > 0$ and all }
\vx\in\bz^3\nz\,.\tag 2.6
$$
\endproclaim

\demo{Proof} By transitivity of the action of the stabilizer $H_B$ of
the form $B$ on the level sets
of $B$ in $\br^3\nz$ and by continuity of $B$ at zero, assertion (2.6) is
equivalent to the norm of $h\vx$ being not less than $\vre$  for some
$\vre > 0$ and all  
$\vx\in\bz^3\nz$ and $h\in H_B$. The latter, in view of Theorem 2.1,
is equivalent to the orbit $H_B\bz^3$ being relatively compact in
$\Omega = SL_3(\br)/SL_3(\bz)$. But $H_B = g^{-1}H_Sg$,
therefore the orbit 
$H_B\bz^3$ is  relatively compact  if and only if so is
the orbit $H_Sg\bz^3$. \qed\enddemo

\proclaim{Corollary 2.8} Let $B$ be a real nondegenerate indefinite
quadratic form in $3$ 
variables. 
If {\rm (2.6)} holds, then 
$
B$ is proportional to  a rational form.
\endproclaim

\demo{Proof} The previous
corollary implies that $H_B\bz^3$ is  compact; but since this orbit
can be identified with $H_B/H_B \cap SL_3(\bz)$, this shows that
$H_B \cap SL_3(\bz)$ is Zariski dense in $H_B$, which is equivalent to
$H_B$ being defined over $\bq$, hence the claim.  \qed\enddemo

To derive Theorem 1.1 from the above corollary, one can then observe that
if $B$ is a real irrational 
nondegenerate indefinite quadratic form in $k$ variables and $l<k$
then $\Bbb R^k$ contains a rational subspace $L$ of dimension $l$
such that the restriction of $B$ to $L$ is irrational nondegenerate
and indefinite (the proof can be found in
\cite{DM1}).  Hence the validity of Theorem 1.1  in the case
$k = 3$ implies the general case. 

\medskip

It is worthwhile to mention that in the paper \cite{Op2} Oppenheim
modified his conjecture  
replacing the claim $\inf_{\vx\in\bz^{k}\nz}|B(\vx)| = 0$ with 
``$0$ is a non-isolated accumulation point of
$B(\bz^k)$'', which he showed to be equivalent to the density of
$B(\bz^k)$ in $\br$. This stronger form of the conjecture was also
proved by Margulis \cite{Ma4}. It is not difficult to derive it from
Theorem 2.5, first reducing  to the case $k = 3$: 
if $B$ is not proportional to  a rational form, the orbit $H_B\bz^3$
is not closed, hence (by non-existence of intermediate subgroups
between  $H_B$ and $SL_3(\br)$) it is dense in $\Omega$, and  the
density of 
$B(\bz^3)$ in $\br$ follows. 

Finally let us briefly mention quantitative extensions of the above
results. For $B$ as above, an open interval $I\subset \br$
and a 
positive $T$ one defines
$$
V_{I,B}(T) \df \{\vx\in\Bbb R^k\mid
B(\vx)\in I, \|\vx\| \le T\}
$$
and
$$
N_{I,B}(T) \df \# \big(\bz^k \cap V_{I,B}(T)\big)\,,
$$
then one has $N_{I,B}(T) \to\infty$ as $T\to \infty$ for every
nonempty $I\subset \br$. The next theorem, a compilation of results
from \cite{DM2} and \cite{EMM}, describes the growth of this counting
function 
comparing it to the volume of $V_{I,B}(T)$. 

 \proclaim{Theorem 2.9} Let  $B$ be a real  nondegenerate indefinite
quadratic form  of
signature $(m,n)$, $m+n > 2$, $m \le n$, which is not proportional
to a 
rational form, and  
let $I\subset \br$ be a  nonempty open interval. 
Then 

{\rm (a) \cite{DM2}} $
\dsize\liminf_{T \to \infty} \frac{N_{I,B}(T)}{ |V_{I,B}(T)|} = 1\,;
$

{\rm (b) \cite{EMM}} if $n \ge 3$, then $\dsize
\lim_{T \to \infty} \frac{N_{I,B}(T)}{ |V_{I,B}(T)|} = 1\,.
$
\endproclaim

In the exceptional cases, i.e.~for forms of signature $(2,1)$ and
$(2,2)$, there are counterexamples showing that $N_{I,B}(T)$ can grow
like const$\cdot |V_{I,B}(T)| (\log T)^{1-\vre}$, and it is
proved in \cite{EMM} that  const$\cdot |V_{I,B}(T)| \log T$ is an
asymptotically exact upper bound for $N_{I,B}(T)$.

For both parts the crucial step is to approximate the counting
function $N_{I,B}(T)$ by values of integrals of certain functions
along orbits in the space of lattices. Part (a)
relies upon Ratner's uniform distribution theorem (a refinement of
Theorem 
2.5 for one-parameter unipotent subgroups), while the second part
involves delicate estimates based on combinatorics of lattices. See \cite{DM1, 
DM2, EMM, Ma7} 
for more 
details on the proofs and further refinements and generalizations.

\head{3. Metric linear \da\ and lattices}\endhead

The basic object to study in this section will be the set $\mr$ of
$m\times n$ real matrices. The word ``metric'' refers to considering
solution sets of \di\ inequalities in terms of the Lebesgue measure, or,
when the sets are of measure zero, for finer analysis, in terms of the \hd. 
We refer the reader to the books  \cite{C, H, S3, Sp3} for
a detailed exposition. 

In what follows, 
$\psi(\cdot)$  will be a positive non-increasing  function $\br_{\sssize
+}\mapsto \br_{\sssize
+}$.  We are going to use it to measure the precision of approximation
of a real number $\alpha$ by rational numbers as follows: we would
like the fractional part of $\alpha q$ to be not bigger than
$\psi(|q|)$ for infinitely many $q\in \bz$ in order to call $\alpha$
"sufficiently well approximable" (this notion being dependent on
$\psi$). More precisely, let us  say that $\alpha$ is 
{\sl $\psi$-approximable} if  there are infinitely many
$q\in \bz$ such that 
$$
 |\alpha q + p|   \le \psi(|q|) \quad \text{for some
}p\in\bz\,. 
$$
In order to consider a matrix analogue of this notion, one needs to
choose a norm on $\br^k$ (we will do it by setting $\|\vx\| =
\max_{1\le i \le k}|x_i|$). Then one says that  a matrix \amr\
(viewed  
as a  system   of $m$ linear forms  in $n$ variables)
is 
{\sl $\psi$-approximable} if  there are infinitely many
$\vq\in \bz^n$ such that 
$$
 \|Y\vq + \vp\|^m   \le \psi(\|\vq\|^n) \quad \text{for some
}\vp\in\bz^m\,.  \tag 3.1
$$
The above normalization (raising norms in the power equal to the
dimension of the space, instead of the traditional $ \|Y\vq + \vp\|
\le \psi(\|\vq\|)$ 
as in \cite{Dod} or \cite{BD}) is convenient for many reasons: in our
opinion it  makes the structure 
more transparent and less dimension-dependent, and simplifies the
connection  with homogeneous flows.

The whole theory starts from a positive result of Dirichlet, namely

\proclaim{Theorem 3.1}  Every \amr\ 
is 
$\psi_0$-approximable, where $\psi_0(x) = \frac1x$. 
\endproclaim  

Clearly the faster $\psi$ decays, the smaller is the set of $\psi$-approximable
matrices. The next theorem, Groshev's \cite{Gr}
generalization of earlier results of Khintchine,
provides the zero-one law 
for the Lebesgue 
measure of this set:

\proclaim{Theorem 3.2 {\rm (The Khintchine-Groshev Theorem)}}
Almost every (resp.~almost no)  
\amr\
is 
$\psi$-approximable,  provided the sum $
{\sum_{l = 1}^\infty {\psi(l)}} 
$
 diverges (resp.~converges). 
\endproclaim

Now say that \amr\ is {\sl \ba\/} if it is not
$c\psi_0$-approximable for some $c > 0$; that is, 
if there exists $c>0$ such that $\|Y\vq + \vp\|^m\|\vq\|^n
\ge c$ for all ${\vp\in\bz^m}$ and
all but finitely many ${\vq\in\bz^n}$ (equivalently: all
$\vq\in\bz^n\nz$). 
 
Note that in the case  $m = n = 1$,  $\alpha\in\br$  is \ba\ if and
only if 
coefficients in the continued fraction expansion 
of $\alpha$ are bounded. Using continued fractions, Jarnik proved in
1928  that \ba\ numbers form a set of \hd\ one; for arbitrary
$m,n$ the 
corresponding fact, i.e.~full \hd\ of \ba\ systems,  was established by
Schmidt in 1969 \cite{S2}.  

The following interpretation of this property in terms of homogeneous
dynamics is due to Dani. Throughout this section we will fix
$m,n\in\bn$ and put $k = \mn$. 

\proclaim{Theorem 3.3 \cite{D1}} \amr\ is
\ba\ iff the 
trajectory \linebreak $\{g_tL_{ Y}\bz^k\mid t\in\br_+\}$, with $\la$ as
in {\rm (1.3)} and $g_t$ as in {\rm (2.1)}, 
 is
bounded in the space $\Omega$ of unimodular lattices in $\br^k$. 
\endproclaim

Instead of giving the proof (which, besides the original paper
\cite{D1} can be found  in
\cite{K2, K3}) let us point out the similarity between the
above theorem 
and Lemma 2.7. Indeed, denote by $S$ the function on $\br^k$ given
by
$$
S(x_1,\dots,x_k) =
\max(|x_1|,\dots,|x_m|)^m\max(|x_{m+1}|,\dots,|x_k|)^n\,.
$$
Then $Y$ is \ba\ iff  for some $c >0$ one has 
$$
S(\la\vx) \ge c\text{ for all }\vx =
(\vp,\vq)\in\bz^m \times (\bz^n\nz)\,.\tag 3.2
$$ 
Furthermore, the one-parameter group $\{g_t\}$ as in
(2.1) is essentially (up to the compact part) the stabilizer of $S$,
and,  as in the proof of Lemma 2.7, one can show 
that (3.2) is equivalent to the norm of $g_t\la\vx$ being bounded away
from zero for all $\vx\in\bz^k\nz$ and $t\ge 0$,
that is, to the statement that
$\{g_tL_{
Y}\bz^k\mid t\in\br_+\}\cap \Omega_\vre = \varnothing$ for some
$\vre > 0$. 

From the above theorem and the aforementioned result of Schmidt,
Dani derived   

\proclaim{Corollary 3.4} The set 
$$
\big\{\Lambda \in \Omega\mid \{g_t\Lambda\mid t\ge 0\} \text{ is bounded} \big\}\,,
$$
with $\{g_t\}$ as in {\rm (2.1)}, has full \hd. \endproclaim

\demo{Proof} Indeed, any
$\Lambda\in\Omega$ can be written as  $\left(\matrix
B & 0  \\
C & D
\endmatrix \right)  L_{
Y} \bz^k$, therefore one has
$$
g_t\Lambda =  g_t\left(\matrix
B & 0  \\
C & D
\endmatrix \right) g_{-t}\cdot g_tL_{ Y} \bz^k\,. 
$$
But as we saw in one of the examples of \S 2, $\{L_{ Y} \mid
Y\in\mr\}$ is the \ehs\ of $G$ relative to $g_1$; thus the conjugation
of the neutral and contracting parts plays no role and the trajectory
$g_t\Lambda$ is bounded iff so is $g_tL_{ Y}\bz^k$. \qed\enddemo

A possibility to generalize the statement of the last corollary (to
actions of other one-parameter groups on other \hs s) was mentioned by
Dani in \cite{D2} and later conjectured by Margulis
\cite{Ma5, \rm
Conjecture (A)}. The latter conjecture was
settled by Margulis and the author in 1996. Let us state here the following 
weakened version:

\proclaim{Theorem 3.5 \cite{KM1, K4}} Let  $G$ be a Lie
group, $\Gamma$ a lattice in $G$, $F = \{g_t\mid t\ge
0\}$ be a one-parameter subsemigroup of $G$ consisting of 
semisimple\footnote{$g\in G$ is called {\sl semisimple} if the operator $\Ad_g$ is diagonalizable over $\bc$}
 elements,
 and let 
$H = H_{\sssize +}(g_1)$ be the expanding horospherical subgroup
corresponding to $F$. Assume in addition that
 the $F$-action  on $\ggm$ is mixing.
Then for any closed $F$-invariant null 
subset $Z$ of $\ggm$ and any $x\in\ggm$, the set
$$
\{h \in H \mid Fhx\text{ is bounded and }\overline{Fhx}\cap
Z = \varnothing\}
$$
 has full \hd. In particular, if $\{g_t\}$ is partially hyperbolic,
then the set 
$
\{x \in \ggm \mid Fx\text{ is bounded and }\overline{Fx}\cap
Z = \varnothing\}
$ has full \hd. 
\endproclaim

Note that abundance of exceptional orbits is a feature of
many {\sl chaotic} dynamical systems. See e.g.~\cite{AN1, AN2,
D4, Dol1, U}.
In the situation of Theorem 3.5, the construction of 
bounded orbits 
(or, more generally, orbits staying away from a fixed part of the
space) comes from uniform distribution of images of expanding leaves
(Theorem 2.3). More precisely, first one reduces the problem to the
case (2.4), and then considers natural ``rectangular''
partitions of $H$  (called {\sl tessellations\/} in \cite{KM1} and
\cite{K1}) and studies their behavior under the automorphism
$h\mapsto g_thg_{-t}$ of $H$. Theorem 2.3 is used to show than one
can cover the set of ``bad'' points by relatively small number of
rectangles. Then those rectangles are used to create a Cantor set
consisting of points with orbits avoiding  $Z$ and staying within a
compact subset of $\ggm$. See \cite{KM1, K1, K4} for details and
generalizations\footnote{In particular, it follows from the methods 
of \cite{KM1} that one can remove the assumption of semisimplicity of 
elements of $F$, but then one needs the $F$-action to be exponentially 
mixing, that is, (2.2) must hold for any $\varphi, \psi$ as in Theorem 2.2.}

\medskip 

So far we have illustrated the impact of ideas coming from \da\ to
ergodic theory. On the other hand, Theorem 3.5 and Dani's
correspondence (Theorem 3.3) can be used as an
alternative proof of the aforementioned result of Schmidt on
abundance of \ba\ systems of linear forms. What follows is another
application to number theory, which produces a new result and
demonstrates the power of ideas relating the two fields. 

Let us consider an inhomogeneous twist of approximation of real
numbers by rationals. Instead of just one real number $\alpha$
take a pair $\langle \alpha,\beta \rangle$, consider an {\sl affine 
form\/} $x\mapsto \alpha x +\beta$ and look at fractional parts of its
values at integers. Similarly, a system  of $m$ affine forms in $n$
variables  will be then given by   
a pair $\aa$, where \amr\ and $\va\in\br^m$. Let us  denote by $\mtr$ the direct product of $\mr$ and $\br^m$.
Now say that a system  of affine forms  given by  
\amtr\ is 
{\sl $\psi$-approximable} if  there are infinitely many
$\vq\in \bz^n$ such that 
$$
 \|Y\vq + \va +\vp\|^m   \le \psi(\|\vq\|^n) \quad \text{for some
}\vp\in\bz^m\,, 
$$
and {\it \ba\/} if  it is not
$c\psi_0$-approximable for some $c > 0$; that is, 
there exists a constant $\tilde c>0$ such that for every ${\vp\in\bz^m}$ and
all but finitely many ${\vq\in\bz^n\nz}$ one has 
$$\|Y\vq + \va + \vp\|^m\|\vq\|^n  > \tilde c \,.
$$ 

It can be proved (and follows from an 
inhomogeneous version of the Khintchine-Groshev Theorem, see \cite{C})
that the set of \ba\  \amtr\ is of measure zero. However, all known
examples of \ba\  \amtr\ belong to a countable union 
of proper submanifolds of $\mtr$, hence  form a set of positive Hausdorff
codimension. Yet a modification of the dynamical approach described
above works in this case as well.  Namely, one considers
a collection of 
vectors 
$$
\left\{\left.\left(\matrix Y\vq + \va + \vp \\  \vq \endmatrix
\right)\right| \vp \in\bz^m,\,\vq \in\bz^n\right\} = L_{
Y}\bz^{k} + \left(\matrix \va \\ 0\endmatrix
\right) \,,
$$
which is an element of the space $\hat \Omega = \hat G/\hat
\Gamma$ of {\sl affine lattices} in $\br^k$, where $$
\hat G\df \text{\rm
Aff}(\br^k) =  G\ltimes\br^k\text{ and }\hat
\Gamma \df \Gamma\ltimes\bz^k\,.
$$
In other words, 
$$
\hat\Omega\cong\{\Lambda +\vw\mid \Lambda\in\Omega,\ \vw\in\br^k\}\,.
$$

Note that the quotient topology on $\hat\Omega$
coincides with the natural topology on the space of affine lattices:
that is, $\Lambda_1 +\vw_1$ and $\Lambda_2 +\vw_2$ are close to each
other if so are $\vw_i$ and  the generating elements of
$\Lambda_i$. Note also that  $\hat \Omega$ is  non-compact and has finite Haar
measure, and that  $\Omega$ (the set of {\sl true} lattices) can be
identified with a subset of $\hat \Omega$ (affine lattices containing
the zero vector). Finally,  $g_t$ as in (2.1) acts on $\hat \Omega$,
and it is not hard to show that the
\ehs\ corresponding to $g_1$ is exactly the set of all elements of
$\hat G$ with linear part $L_{
Y}$ and translation part $\left(\matrix \va \\ 0\endmatrix
\right) $, \amr\ and $\va\in\br^m$.

Now, for $\vre > 0$, define 
$$
\hat \Omega_\vre \df\big\{\Lambda \in \hat\Omega \bigm|  \|\vx\| < \vre\text{ for some
}\vx\in\Lambda\big\}\,.
$$
Then $\hat \Omega\ssm\hat \Omega_\vre$ is a closed  (non-compact) set disjoint from
$\Omega$.  

\proclaim{Theorem 3.6 \cite{K4}}  Let $F = \{g_t\mid t\ge 0\}$
be as in {\rm 
(2.1)}. Then 
$$ 
F\left(L_{
Y}\bz^{k} + \left(\matrix \va \\ 0\endmatrix
\right)\right) \quad \text{is bounded and stays away from }\Omega
$$
$$
\Downarrow
$$
$$
F\left(L_{
Y}\bz^{k} + \left(\matrix \va \\ 0\endmatrix
\right)\right)\subset \hat \Omega\ssm\hat \Omega_\vre\text{ for some }\vre > 0
$$
$$
\Downarrow
$$
$$
\aa\text{ is \ba}
$$
\endproclaim

The proof is basically a slight modification of ideas involved in the
proof of Theorem 3.3. 
It follows from the results of \cite{BM} (see also \cite{Ma5})
that the $F$-action on $\hat \Omega$ is mixing. 
Since $\Omega\subset\hat \Omega$ is closed, null
and $g_t$-invariant, Theorem 3.5 applies and one gets

\proclaim{Corollary 3.7}  The set   of \ba\ \amtr\ has full \hd.
\endproclaim

See \cite{K4} for details, remarks and extensions. 

We close the section by stating a
theorem generalizing Dani's correspondence (Theorem 3.3) to
$\psi$-approximable systems. First we need a simple ``change of
variables'' lemma.

\proclaim{Lemma 3.8}  Fix  $m,n\in\bn$  and $x_0 > 0$, and let
$\psi:[x_0,\infty)\mapsto (0,\infty)$ be a non-increasing continuous
function. Then there exists a unique 
continuous  function  $\vre:[t_0,\infty)\mapsto (0,\infty)$, where
$e^{kt_0} =  x_0^m /\psi(x_0)^n$, such that  
$$
\text{the function} \quad {t\mapsto e^t \vre(t)^n}  \quad \text{is strictly increasing and unbounded}\,,
\tag 3.3a
$$
$$
\text{the function} \quad {t\mapsto e^{-t} \vre(t)^m}  \quad \text{is nonincreasing}\,,
\tag 3.3b
$$
and
$$
\psi\big(e^t \vre(t)^n\big) =  e^{-t} \vre(t)^m\quad\forall\,t\ge
t_0\,.\tag 3.4 
$$
Conversely, given $t_0 \in\br$ and a continuous  function
$\vre:[t_0,\infty)\mapsto(0,\infty)$ such that {\rm (3.3ab)} hold, there exists
a unique continuous non-increasing function
$\psi:[x_0,\infty)\mapsto(0,\infty)$, with $x_0 = e^{t_0} \vre(t_0)^n$,
satisfying {\rm (3.4)}. \endproclaim

See \cite{KM3} for the proof. In many cases one can explicitly solve 
(3.4) to express $\vre(\cdot)$ knowing $\psi(\cdot)$ and vice versa. For
example if $\psi(x) = c\psi_0(x) = c/x$, the equation (3.4) gives  
$ce^{-t} \vre(t)^{-n} = e^{-t} \vre(t)^m$, and one sees that the
corresponding function $\vre$ is constant (more precisely,
$\vre(t) \equiv 
c^{1/{k}}$). Or one can take $\psi(x) =
c\psi_\beta(x)$, where 
$\psi_\beta(x)
= \frac1{x^{1 + \beta}}$, $\beta > 0$; then $\vre(t)$ decreases
exponentially, namely
$$
\vre(t) = c^{\gamma/\beta} e^{-\gamma t}\,, \text{ where } \gamma =
\frac\beta{{(1+\beta)n + m}}\,.
$$

Now we can state a generalization of Theorem 3.3:

\proclaim{Theorem 3.9}  \amr\ is $\psi$-approximable iff
there exist arbitrarily large 
positive  $t$ 
such that 
$
g_{t}L_{ Y}\bz^k \in \Omega_{\vre(t)} \,,
$ 
where   $\{g_t\}$ is as in {\rm (2.1)},
$\la$ as in  {\rm (1.3)}, and $\vre(\cdot)$ is the function
corresponding to $\psi$ as in the previous lemma. \endproclaim  

Loosely speaking, good rational approximations for $Y$ correspond to
far excursions of the orbit into the ``cusp neighborhoods''
$\Omega_{\vre}$. In other words, one can measure the ``growth rate''
of the orbit in terms of hitting the sets $\Omega_{\vre(t)}$ in time
$t$ for infinitely many $t\in \bn$, and fast-growing orbits would
correspond to systems approximable with a fast-decaying approximation
function. 

It is shown in \cite{KM3} how the above
correspondence provides an alternative (dynamical) proof of Theorem
3.2. More precisely, one can use ergodic properties of the
$g_t$-action on $\Omega$ (exponential decay of correlations, see
Theorem 2.2) to prove the following

\proclaim{Theorem 3.10} Let $\vre(\cdot)$ be any positive
function. Then for almost all (resp.~almost no) 
$\Lambda\in\Omega$ 
  one has  $g_t \Lambda \in \Omega_{\vre(t)}$ for infinitely
many  $t\in\bn$,  provided the sum
$$
\sum_{t = 1}^\infty \vre(t)^k \tag 3.5
$$
 diverges (resp.~converges). \endproclaim

We remark that the ratio ${\mu(\Omega_{\vre})}/{\vre^k}$ is
shown in \cite{KM3} to be bounded from both sides; therefore the sum
(3.5) is finite/infinite iff so is $\sum_{t = 1}^\infty
\mu(\Omega_{\vre(t)})$. This places the above theorem in the rank of
Borel-Cantelli type results. See \cite{KM3} for
generalizations and applications, and \cite{CK, CR, Dol2, Ph,
Su} for other 
results of similar flavor.

Another application of the correspondence of Theorem 3.9 will be given in
the next section.

\head{4. \da\ on manifolds}\endhead

We start from the setting of the previous section but specialize to
the case $m = 1$; that is, to \da\ of just one linear form  given by
$\vy\in\br^n$.  Recall that Theorem 3.2 says  
that whenever $\sum_{l = 1}^\infty {\psi(l)}$ is finite, almost  every 
$\vy$ 
is not
$\psi$-approximable; that is, the inequality 
$$
|{\vq}\cdot{\vy} + p|
\leq\psi(\|{\vq}\|^n)\tag 4.1
$$
 has at most finitely many solutions. It is
instructive to sketch an elementary proof:  for fixed $p,\vq$, the set
of $\vy$ satisfying
(4.1) is a $\frac{\psi(\|{\vq}\|^n)}{\|{\vq}\|}$-neighborhood of a
hyperplane 
$$
{\vq}\cdot{\vy} + p = 0\,;\tag 4.2
$$
 thus if one restricts $\vy$ to lie
in  $[0,1]^n$ (or any other bounded subset of $\br^n$), the set of
solutions will  have measure at most const$\cdot
\frac{\psi(\|{\vq}\|^n)}{\|{\vq}\|}$. Since there are at most const$\cdot
\|{\vq}\|$ admissible values of $p$,  the sum of measures of all
sets of solutions is at most 
$$
\sum_{\vq\in\bz^n}\psi(\|{\vq}\|^n) \asymp \sum_{l = 1}^\infty l^{n-1}
\psi(l^n) \asymp\sum_{l = 1}^\infty \psi(l)\,,
$$
and the proof is finished by an application of the Borel-Cantelli
Lemma.

 Recall that $\psi_\beta(x) = x^{-(1+\beta)}$, $\beta > 0$, was one of
the examples of functions realizing the convergence case in the
Khintchine-Groshev Theorem. Say that $\vy\in\br^n$ is {\sl very well
approximable (VWA)\/} if it is 
$\psi_\beta$-approximable for some   $\beta > 0$.  Thus almost all
${\vy}\in\br^n$ are not  VWA.

Now consider the following problem, raised by Mahler in
1932 \cite{M}: is it
true that for almost all $x\in\br$ the inequality
$$
|p + q_1x + q_2 x^2 + \cdot + q_nx^n|
\leq \|{\vq}\|^{-n(1+\beta)}
$$
 has at most finitely many solutions? In other words, for
a.e.~$x\in\br$, the $n$-tuple 
$$
\vy(x) = (x,x^2,\dots,x^n)\tag 4.3
$$ is not VWA. 
The proof presented above does not work, since this time one has to
estimate the measure of intersection of the curve (4.3) with the sets
of solutions of inequalities (4.1), and for some choices of $p,\vq$
(namely for those which make the hyperplane (4.2) nearly tangent to
the curve (4.3)) it is hard to produce a reasonable estimate. 

This problem stood open for more than 30 years until it was solved in 1964
by \sp\ \cite{Sp1, Sp2}. Earlier several special cases were
considered, and, quoting 
\sp's survey paper, the problem  rapidly revealed itself to be
non-trivial and involving 
``deep and complicated phenomena in which arithmetical properties of
numbers are closely entangled with combinatorial-topological
properties''  of the curve. The solution to Mahler's problem  has
eventually led to the development of a new branch 
of metric number theory, usually referred to as ``\da\ with dependent
quantities'' or ``\da\ on manifolds''. 
We invite the reader to look at \sp's monographs \cite{Sp2,
Sp3} and a recent book \cite{BD} for a systematic exposition
of the field. 

\medskip 

Mahler's problem and its generalizations have several motivations. The
original motivation of Mahler comes from transcendental number
theory. Indeed, the  $n$-tuple 
(4.3) is not VWA if and only if for every $\beta > 0$ there are at
most finitely many polynomials $P\in\bz[x]$ with degree at most $n$
such that $|P(x)| < h(P)^{-n(1+\beta)}$, where $h(P)$ is the height of
$P$; loosely speaking, $x$ is ``not very algebraic'', and  the
affirmative  solution to the problem shows that almost all $x$ are
such. 

Another motivation comes from KAM theory: it is known that behavior of 
perturbation of solutions of ODEs is related to \di\ properties of
coefficients. If the latter are restricted to lie on a curve or
submanifold of $\br^n$, it may be important to know that almost all
values have certain
approximation properties. See
\cite{\rm de la Llave's
lectures, this volume} and \cite{BD, \rm Chapter 7}. 

However, from the author's personal viewpoint, the appeal of this
branch of number theory lies in  its existing and potential
generalizations. In a sense, the affirmative solution to Mahler's
problem shows that a certain
property of $\vy\in\br^n$ (being not VWA) which holds for generic
$\vy\in\br^n$  in fact holds for generic points on the curve (4.3). 
In other words, the curve inherits the above \di\ property from the
ambient space, unlike, for example, a line $\vy(x) = (x,\dots,x)$ --
it is clear that every point on this line is VWA. This gives rise to
studying other subsets of $\br^n$ and other \di\ properties, and
looking at whether this inheritance phenomenon takes place.

\medskip 

Note that the curve
(4.3) is not contained in any affine subspace of $\br^n$ (in other
words, constitutes an essentially $n$-dimensional object). The latter
property, or, more precisely, its infinitesimal analogue, is
formalized in the following way. Let $V$ be an open subset of
$\br^d$. Say that an $n$-tuple  $\vf = (f_1,\dots,f_n)$ of  $C^l$ 
functions $V\mapsto \br$ is {\sl nondegenerate at $x\in V$\/} if 
the space $\br^n$ is spanned by partial derivatives of $\vf$ at
$x$ of order up to $l$. 
If $M\subset \br^n$ is a $d$-dimensional smooth submanifold, we will say
that $M$ is {\sl nondegenerate at $\vy\in M$} if any (equivalently,
some)  diffeomorphism
$\vf$ between an open subset $V$  of $\br^d$ and a neighborhood of
$\vy$ in $M$ is   nondegenerate at $\vf^{-1}(\vy)$. We will say that
$\vf:V\to \br^n$ (resp.~$M\subset \br^n$) is {\sl nondegenerate\/} 
if it is nondegenerate at almost every point of $V$  (resp.~$M$, in
the sense of the 
natural measure class on $M$). 
If the functions $f_i$ are analytic, it is easy to see that the linear
independence of $1,f_1,\dots,f_n$  over $\br$ in  $V$ is
equivalent to all points of $M = \vf(U)$ being {nondegenerate}.  Thus
the above nondegeneracy  condition can be viewed as an infinitesimal
version of not lying in any proper affine subspace of $\br^n$. 

It appears that many known and anticipated results in the field fall
in the framework of the following
vague

\proclaim{Meta-Conjecture} ``Any'' \di\ property of vectors in an ambient space
(e.g.\ $\br^n$) {which} {holds for almost}
{all points in this space should 
hold for generic
points on a nondegenerate smooth submanifold $M$ of} {the
space}. \endproclaim 

It was conjectured in 1980 by
\sp\ \cite{Sp4, \rm Conjecture
H$_1$} that almost all points on a
nondegenerate analytic submanifold of $\br^n$ are not VWA. This
conjecture was supported 
before and after 1980 by a number of partial results, one of the first
being Schmidt's proof \cite{S1} for nondegenerate planar curves. 
The general case was settled in 1998 by
Margulis and the author using the dynamical approach. Namely,
the following was 
proved:

\proclaim{Theorem 4.1 {\rm \cite{KM2}}} Let $M$
be a nondegenerate 
smooth submanifold of $\br^n$. Then almost all points of $M$  
are not VWA.
\endproclaim

This is the result we will focus on later in this section. 

\medskip

In another direction, \sp's solution to Mahler's problem was improved
in 1964 by Baker \cite{B1} and later (1984) by Bernik
\cite{Bern, BD}; the latter proved
that  whenever $\sum_{l = 1}^\infty {\psi(l)}$ is finite,  almost all
points of the curve (4.3) are not $\psi$-approximable. And several
years ago Beresnevich \cite{Bere1} proved the divergence
counterpart, thus 
establishing a complete analogue of the
Khintchine-Groshev Theorem for the curve (4.3). 

It turned out that a modification of the methods from \cite{KM2} allows
one to prove the convergence part of the Khintchine-Groshev Theorem
for any nondegenerate manifold. In other words, the following is
true:  

\proclaim{Theorem 4.2} Let $M$ be a nondegenerate
smooth submanifold of $\br^n$ and let $\psi$ be such that
$\sum_{l = 1}^\infty {\psi(l)}$ is finite. Then almost all
points of $M$   
are not $\psi$-approximable.
\endproclaim

This is proved in \cite{BKM} and also independently in
\cite{Bere2}. A work on the divergence case is currently in
progress. (We
note that the main result of \cite{BKM} has a stronger
``multiplicative'' 
version (see Theorem 5.3) which is currently not doable by classical \sp-style methods
developed 
in \cite{Bere2}.)

In this survey we will indicate a proof of Theorem 4.1 by first
restating it in the language of 
flows on the space of lattices. For this we set $k = n+1$ and
look at the 
one-parameter group 
$$
g_t = \text{\rm
diag}(e^{t},e^{-t/n},\dots,e^{-t/n})\,.\tag 4.4
$$
acting on $\Omega = SL_{k}(\br)/SL_{k}(\bz)$, and given $\vy\in\br^n$,
consider 
$
L_{\vy} \df \left(\matrix
1 & \vy^{\sssize T} \\
 0 & I_n
\endmatrix \right)
$
(cf.~(2.1) and (1.3)). It follows from Theorem 3.9 and the example
discussed afterwards that $\vy\in\br^n$ is VWA
iff for some $\gamma > 0$   there exist arbitrarily large
positive  $t$ 
such that 
$$
g_{t}L_{\vy}\bz^{k} \in \Omega_{e^{-\gamma t}} \,.\tag 4.5
$$ 
Equivalently, for some $\gamma > 0$   there are infinitely many
$t\in\bn$ such that (4.5) holds.

With this in mind, let us turn to the setting of Theorem 4.1. 
Namely let $V$ be an open subset of $\br^d$ and  $\vf = (f_1,\dots,f_n)$  an $n$-tuple  of  $C^k$ 
functions $V\mapsto \br$ which is nondegenerate at almost every point
of $V$.  The theorem would be
proved if we show that  for any $\gamma > 0$ the set
$$
\{x\in V\mid g_{t}L_{\vf(x)}\bz^{k} \in \Omega_{e^{-\gamma
t}}\text{ for infinitely many }
t\in\bn\} 
$$
has measure zero. In other words, a submanifold $\vf(V)$ of $\br^n$
gives rise to a submanifold $ L_{\vf(V)}\bz^{k}$ of the space of
lattices, and one needs to show that the growth rate of generic orbits
originating from this submanifold is consistent with the growth rate of
an orbit of a generic point of $\Omega$ (see Theorem 3.10 for an
explanation of why lattices $\Lambda$ such that $g_{t}\Lambda \in
\Omega_{e^{-\gamma t}}\text{ for infinitely many }
t\in\bn $ form a null subset of $\Omega$).

Now one can use the Borel-Cantelli Lemma to reduce Theorem 4.1 to
the following statement:

\proclaim{Theorem 4.3} Let $V$ be an open subset of $\br^d$ and  $\vf
= (f_1,\dots,f_n)$  an $n$-tuple  of  $C^k$  
functions $V\mapsto \br$ which is nondegenerate at $x_0\subset V$.
Then there exists a neighborhood $B$ of $x_0$ contained in $V$ such
that  for any $\gamma > 0$ one has
$$
\sum_{t = 1}^\infty |\{x\in B\mid g_{t}L_{\vf(x)}\bz^{k} \in \Omega_{e^{-\gamma
t}}\} | < \infty\,.\tag 4.6
$$
\endproclaim

Here is the turning point of the argument: $t$-dynamics gives way to
$x$-dynamics, namely, a natural way to demonstrate (4.6) is to fix
$t$ and think of the set $\{g_{t}L_{\vf(x)}\bz^{k}\}$ as of an
orbit of certain action (not a group action!), the goal being to prove
that a substantial part of this ``orbit'' lies outside of  ``cusp
neighborhoods'' $\Omega_{e^{-\gamma
t}}$ uniformly for all $t$. What immediately comes to mind is the
recurrence property of unipotent orbits, that is, Theorem 2.4. 
And it turns out that a modification of the argument used to prove the
latter theorem allows one to 
estimate the amount of ``time'' $x$ that the ``trajectory''
$x\mapsto g_{t}L_{\vf(x)}\bz^{k}$ spends ``close to
infinity'' in $\Omega$. More precisely, the following can be proved:

\proclaim{Theorem 4.4} Let $V$, $\vf$  and $x_0$ be as in Theorem 4.3.
Then there exists a neighborhood $B$ of $x_0$ contained in $V$ and
constants $C,\rho > 0$ such
that  for any $\vre > 0$ and any positive $t$ one has
$$
|\{x\in B\mid g_{t}L_{\vf(x)}\bz^{k} \in \Omega_{\vre}\} | <
C \vre ^{1/dn}  |B|\,.\tag 4.7
$$ 
\endproclaim

It is straightforward to verify that Theorem 4.3 follows from the
above uniform estimate. 

It remains to explain why the behavior of the curve
$x\mapsto g_{t}L_{\vf(x)}\bz^{k}$ is similar to that of the
unipotent orbit. In fact it has been understood a long time ago that
the main property of the unipotent actions on which the recurrence
estimates are based is the polynomial
dependence of   $u_x$  on $x$ (see \cite{Ma2, Sh}). One may wonder 
what 
is so special about polynomials -- and it turns out that the crucial
property is roughly ``not making very sharp turns''. More precisely,
here is the definition motivated by the analysis of the proofs in
\cite{Ma2} (see also \cite{EMS}): for $C,\alpha > 0$ say that a
continuous function $f$ on an 
open set $V\subset 
\br^d$ is {\it \cag\ on $V$\/} if for any open ball $B\subset V$
and any positive $\vre$ one has 
$$
\big|\{x\in B\bigm| |f(x)| < \vre\cdot{\sup_{x\in B}|f(x)|}\}\big| \le
C\vre^\alpha |B|\,. 
$$

In other words, a good function which takes small values on a big part
of a ball is not allowed to grow very fast on the remaining part
of the ball. The main example is provided by polynomials:

 \proclaim{Lemma 4.5} Any
polynomial $f\in\br[x]$ of degree not greater than $k$ is  
$(4k, 1/k)$-good on $\br$.
\endproclaim
 
This easily follows from Lagrange's
interpolation formula, see \cite{DM2, KM2}. The next theorem (the main
result of \cite{KM2}) 
therefore provides a 
generalization of Theorem 2.4 to polynomial trajectories on $\Omega$. 
To state it we need to introduce some notation.   If $\Delta$ is a
discrete subgroup  of $\br^k$ (not necessarily a lattice) generated by
$\vv_1,\dots,\vv_l$, let us measure its norm, $\|\Delta\|$, by the norm
of the exterior product $\vv_1\wedge\dots\wedge\vv_l$. For this one
needs to extend the norm from $\br^k$ to its exterior algebra. If
$\ve_1,\dots,\ve_k$ are standard  base vectors of $\br^k$,  the
elements $\ve_{\sssize I} \df \ve_{i_1}\wedge\dots\wedge \ve_{i_l}$,
$I = \{i_1,\dots,i_l\} \subset
\{1,\dots,k\}$ form a basis of $\bigwedge^l(\br^k)$. Since 
\di\ applications call for the supremum norm, we will extend
$\|\cdot\|$ to $\bigwedge(\br^k)$ by setting  $\|\sum_Iw_I
\ve_I\| = \max_I|w_I|$. 

Now let us consider a  curve in $\Omega$ given by $x\mapsto
h(x)\bz^k $, where $h$ is some function from $\br^d$ to
$GL_k(\br)$. It turns out that in order to understand its recurrence
properties  one has to keep an eye on norms
of all discrete subgroups of $h(x)\bz^k $; in particular, it will be
necessary to prove that all those norms (as functions of
$x$) are \cag\ for some $C,\alpha$. In fact, it will suffice to look
at the coordinates of $h(x)(\vv_{1}\wedge\dots\wedge \vv_{l})$ where
$\vv_{1},\dots, \vv_{l}$ form a basis of $\Delta\subset \bz^k$; one
can easily show that if all components of a vector function are \cag,
the norm of this  function is also \cag.

 \proclaim{Theorem 4.6} Let $d,k\in\bn$, $C,\alpha > 0$, $0 < \rho  \le 1/k$,  and let a ball $B = B(x_0,r_0)\subset \br^d$ and a map $h:\tilde B \to GL_k(\br)$ be given, where $\tilde B$ stands for $B\big(x_0,3^kr_0\big)$.  
Assume that for any subgroup $\Delta$ of $\bz^k$, 
\roster
\item"(i)" the function $x\mapsto \|h(x)\Delta\|$ is \cag\ on $\tilde B$;
\item"(ii)" $\sup_{x\in B}\|h(x)\Delta\| \ge \rho $.
\endroster
Then 
 for any  positive $ \vre \le \rho$ one has
$$
\left|\{x\in B\mid h(x)\bz^k \in\Omega_\vre\}\right| \le \text{\rm
const}(d,k)\cdot \left(\frac\vre \rho \right)^\alpha  |B|\,.\tag 4.8
$$
\endproclaim

 \proclaim{Corollary 4.7}  For any lattice  $\Lambda$ in $\br^k$ there exists a constant $\rho = \rho(\Lambda) > 0$ such that for any one-parameter unipotent subgroup $\{u_x\}_{x\in\br}$   of $SL_k(\br)$, for any $T > 0$ and any $\vre \le \rho$, one has
$$
\left|\{0 < x < T\mid u_x\Lambda\in\Omega_\vre \}\right| \le
\text{\rm const}(k) \left(\frac\vre \rho \right)^{1/k^2}T\,.\tag
4.9
$$
\endproclaim

This is clearly a quantitative strengthening of Theorem 2.4, with an
explicit estimate of $\delta$ in terms of $\vre$.

\demo{Proof} Write $\Lambda$ in the form  $g\bz^k$ with $g\in
GL_k(\br)$, and denote by $h$ the function $h(x) = u_xg$. For any
$\Delta\subset \bz^k$ with basis $\vv_{1},\dots, \vv_{j}$,
the coordinates of $h(x)(\vv_{1}\wedge\dots\wedge \vv_{j})$
will be polynomials in $x$ of degree not exceeding $k^2$. Hence  the functions
$x\mapsto\|h(x)\Delta\|$ will be
$(C,{1/k^2}\big)$-good on $\br$, where $C$ is a constant depending only
on $k$. Now let $\rho \df \min\big(1/k,\inf_{\Delta \in \Cal
L(\bz^k)}\|g\Delta\|\big)$, positive by the discreteness of $\Lambda$
in $\br^k$. Then $\|h(0)\Delta\| \ge \rho$ for any
$\Delta\subset \bz^k$, therefore, with the the substitutions $B = (0,T)$,
$\alpha = {1/k^2}$ and $d = 1$  assumptions (i) and (ii) of Theorem
4.6 are satisfied, and one immediately gets (4.9) from
(4.8). \qed\enddemo 
 
As was mentioned in \S 2, the proof of Theorem 4.6 is based on
delicate combinatorial (partially ordered) structure of the space of 
lattices, and the reader is referred  to \cite{KM2} or \cite{BKM} (most
of  the
ideas are borrowed from \cite{Ma2} and \cite{D3}). Assuming the latter
theorem, we  
conclude by presenting a 

\demo{Sketch of proof of Theorem 4.4}  Take a positive $t$ and
consider $h(x) \df g_{t}L_{\vf(x)}$; clearly all one needs to
prove (4.7) is to check conditions (i) and (ii) of Theorem 4.6 for
every $\Delta\subset \bz^k$. An
elementary 
computation shows
that the coordinates of $h(x)(\vv_{1}\wedge\dots\wedge \vv_{l})$ for
any choice of vectors $\vv_i$ are linear combinations of functions
$f_1,\dots,f_n$ and $1$. Consider first the case of original Mahler's
conjecture, with $d = 1$ and $f_i(x) = x^i$. Then, as in the proof of
Corollary 4.7, condition (i) is automatic due to Lemma 4.5. Further,
 a straightforward computation of the action of  $h(x) $ on
exterior products of vectors in $\br^{k}$ shows that at least one
 coefficient of at least one polynomial arising as a coordinate must
have absolute value not less than 
$1$. This implies that for every
interval $B$ there exists a constant $\rho$ (independent of $t$)
such that (ii) holds. 

It remains to pass from this special case to the general situation of
functions $f_1,\dots,f_n$ on $\br^d$ coordinatizing a nondegenerate
submanifold of $\br^n$. Here one basically has to show that locally
these functions behave like polynomials. Indeed, the following was
proved in \cite{KM2}:

\proclaim{Lemma 4.8} Let $\vf = (f_1,\dots,f_n)$ be a $C^l$ map from
an open subset $V$ of $\br^d$ to $\br^n$, and let $x_0\in V$ be such
that $\br^n$ is spanned by partial derivatives of $\vf$ at $x_0$ of
order up to $l$. Then there exists a neighborhood $U\subset V$ of
$x_0$ and positive $C$  such that any linear combination of $
1,f_1,\dots,f_n$ is $(C,1/dl)$-good on $V$. \endproclaim 

Now to finish the proof one simply has to choose $U$ according to the
above lemma, then pick a ball $\tilde B$ centered at $x_0$ and
contained in $U$, and finally take $B$ to be a concentric ball with
radius $3^{k}$ times smaller. This implies condition (i) with
$\alpha = 
1/dn$ and some constant $C$ independent of $\Delta$ and $t$,
and (ii) follows as a result of a computation described above:
one shows that  at least one
coordinate  of $h(x)(\vv_{1}\wedge\dots\wedge \vv_{l})$ must have
the form $c_0 + \sum_{i = 1}^nc_if_i(x)$ with
$\max(|c_1|,\dots,|c_n|)\ge 1$, and 
therefore one gets a lower bound (again independent of $t$ and
$\Delta$) for 
$\sup_{x\in B}\|h(x)\Delta\|$. \qed\enddemo

\head{5. Multiplicative approximation}\endhead
We have already seen
in Conjecture 1.2 how the 
magnitude of the integer vector $\vq$ was measured by taking the {\sl
product\/} of 
coordinates rather than the maximal coordinate (that is the norm of
the vector). Let us formalize it by saying, for $\psi$ as before, that
\amr\ is 
{\sl $\psi$-multiplicatively approximable\/} ($\psi$-MA) if  there are
infinitely many 
$\vq\in \bz^n$ such that 
$$
 \Pi(Y\vq + \vp)  \le \psi\big(\Pi_{\sssize +}(\vq)\big) \quad \text{for some
}\vp\in\bz^m\,. 
$$ 
where for $\vx = (x_1,\dots,x_k)\in\br^k$ one defines
$\Pi(\vx) = \prod_{i = 1}^k |x_i|$ and $\Pi_{\sssize +}(\vx) =
\prod_{i = 1}^k \max(|x_i|, 1)$. Clearly any $\psi$-approximable
system of linear forms is automatically $\psi$-MA, but the converse is
not necessarily true. Similarly to the standard setting, one can define
{\sl badly multiplicatively approximable (BMA)\/} and  {\sl very well
multiplicatively approximable (VWMA)\/} systems. It can be easily shown
that almost no \amr\
  are
$\psi$-MA if  the sum 
$$
\sum_{l = 1}^\infty (\log l)^{k - 2}{\psi(l)} 
$$ 
converges (here we again set $k = \mn$); in particular, VWMA
systems form a set of measure zero. 
The converse (i.e.\ a multiplicative analogue of Theorem 1.2) can be
proved using methods of Schmidt; the case $n = 1$ 
is contained in \cite{G}. 

On the other hand, saying that a vector
$\vy\in\br^n$ (viewed as a linear form
$\vq\mapsto \vy\cdot\vq$) is not BMA is equivalent to (1.1);
in other words, Conjecture 1.2 states that no $\vy\in\br^n$, $n \ge
2$, is badly multiplicatively approximable.  
A more general statement that  {no} \amr\ is BMA unless $m = n = 1$ in
fact reduces to 
this conjecture; moreover, as it is the case with Theorem 1.1, it is
enough to prove Conjecture 1.2 for $n = 2$.

It seems natural to bring lattices into the game. In fact, one can
observe 
the similarity between the statements of the Oppenheim and
Littlewoods's conjectures. Indeed, to say that \amr\ is not BMA 
amounts
to saying that $0$ is the infimum of absolute values of a
certain homogeneous polynomial at 
integer points. Similarly to what was done for quadratic
forms, a linear change of variables transforms this 
polynomial into the product of coordinates $\Pi(\vx)$, $\vx\in\br^k$,
and according 
to the scheme developed in the preceding sections, the dynamical
system reflecting \di\ properties of $A$ must come from the group
stabilizing $\Pi$, that is, the full diagonal subgroup $D$ of
$SL_k(\br)$.  
Thus the problems rooted in multiplicative \da\ bring us to {\sl
higher rank actions\/} on $\Omega$ (in implicit form this was already
noticed in the paper \cite{CS} of Cassels and Swinnerton-Dyer). More
precisely, one can state a multiplicative version of generalized
Dani's correspondence (Theorem 3.8), relating multiplicative \di\
properties of $A$  
to orbits of the form $\{g\la\bz^{k}\mid g\in D_+\}$ where $D_+$ is
a certain open chamber in the group $D$. (For a version of
such a correspondence see \cite{KM3, \rm Theorem 9.2}.)

We illustrate this principle by two examples below, where for the sake
of simplicity of exposition we specialize to the case 
$m = 1$ (one linear form
$\vq\mapsto \vy\cdot\vq$, $\vy\in\br^n$), setting $k = n+1$. We will need the following
notation: for $\vt = (t_1,\dots,t_n)\in\br^n$ let us denote $\sum_{i
= 1}^nt_i$ by $t$ and define
$$
g_\vt = \text{\rm diag}(e^{t},e^{-t_1},\dots,e^{-t_n})\in
SL_{k}(\br)\,. \tag 5.1
$$

\subhead Littlewoods's conjecture\endsubhead  One can show
that $\vy\in\br^n$ satisfies (1.1) (that is, it is BMA) iff the
trajectory 
$\{g_\vt\ly\bz^{k}\mid \vt \in\br_+^n\}$  is bounded in the space of
lattices in $\br^{k}$. Thus Conjecture 1.2 is equivalent to the
statement that every trajectory as above is unbounded. In fact, an
argument rooted in the ``Isolation Theorem'' of Cassels and
Swinnerton-Dyer \cite{CS, LW, Ma4} shows that the latter
statement can be 
reduced to the 
following

\proclaim{Conjecture 5.1}  Let $D$ be the subgroup  of diagonal
matrices in
$SL_{k}(\br)$, $k \ge 3$.   Then any relatively compact orbit
$D\Lambda$, $\Lambda$ a lattice in 
$\br^{k}$, is compact. 
\endproclaim

Notice that Theorem 3.5 shows that the above statement does not hold
if $k = 2$. This highlights the difference between rank-one and higher
rank dynamics. 
Note also the similarity between Corollary 2.6 and Conjecture 5.1,
showing that higher rank hyperbolic actions share some features with
unipotent dynamics. In fact, Conjecture 5.1 is a special case of a
more general hypothesis, see
\cite{Ma8, {\rm Conjecture 1}}, which,
roughly speaking, says that for a connected Lie group $G$, a lattice
$\Gamma\subset G$ and a closed subgroup $H$ of $G$, any orbit closure
$\overline{Hx}$, $x\in\ggm$, is an orbit of an intermediate subgroup
$L\supset H$ of $G$ unless ``it has a good reason not to'' (the latter
reasons must be coming from certain one-parameter quotient
actions). See \cite{Ma8, {\rm \S 1}} for more detail.

\subhead Multiplicative approximation on manifolds \endsubhead
Since every VWA vector is VWMA (that is,
$\psi_\vre$-MA for some $\vre > 0$) but not other way around, it is a
more difficult problem to prove that a generic point on a
nondegenerate manifold is not very well multiplicatively
approximable. This has been known as Conjecture 
H$_2$ of \sp\ \cite{Sp4}; the polynomial special case (that is, a
multiplicative strengthening of Mahler's problem) was conjectured by
Baker  in \cite{B2}; both conjectures stood open, except for
low-dimensional special cases, until \cite{KM2} where the following
was proved:

\proclaim{Theorem 5.2} Let $M$ be a nondegenerate
smooth submanifold of $\br^n$. Then almost all points of $M$  
are not VWMA.
\endproclaim

The strategy of the proof of Theorem 4.1 applies with minor changes. 
One shows (see \cite{KM2, \rm Lemma 2.1 and Corollary 2.2} for a
partial result)  that  $\vy\in\br^n$ is VWMA  iff  for some $\gamma > 0$   there  are infinitely many
$\vt\in\bz_+^n$ such that 
$$
g_{\vt}L_{\vy}\bz^{k} \in \Omega_{e^{-\gamma t}}\tag 5.2
$$ 
(here, as defined above, $t$ stands for $\sum_{i
= 1}^nt_i$). Therefore it is enough to use Theorem 4.6  to prove a
modification of the measure estimate of Theorem 4.4 with $g_t$ as in
(4.5)  
replaced by $g_{\vt}$ as in (5.1). 

Finally let us mention a multiplicative version of Theorem 4.2,
proved in \cite{BKM} by a modification of the method described above: 
 
\proclaim{Theorem 5.3} Let $M$ be a nondegenerate smooth 
submanifold of $\br^n$ and let $\psi$ be such that $\sum_{l =
1}^\infty (\log l)^{n - 1}{\psi(l)}$ is finite. Then almost all points
of $M$  
are not $\psi$-multiplicatively
approximable.
\endproclaim

\head Acknowledgements\endhead This survey is based on a
minicourse of lectures given at the AMS Summer Research Institute in
Smooth Ergodic Theory and  
applications (Seattle, 1999). The author is grateful to the organizers
and participants of the workshop, and especially to Alex Eskin for
sharing the responsibility and the fun of giving the course. Thanks
are also due to participants/organizers of the workshop on Ergodic
Theory,  
Rigidity and Number Theory (Cambridge UK, January 2000) where the
presentation of material was tested one more time in a series of
lectures, to Gregory Margulis,
Alexander Starkov and Barak Weiss for helpful discussions, and
to Stella for her endless patience and support.

\enddocument